\documentclass{amsart}

\usepackage{amsxtra,amssymb,multicol,graphicx}
\usepackage{AKstyle}

\usepackage{comment, url}

\usepackage{caption}
\usepackage[labelfont=rm]{subcaption}

\usepackage[pagebackref=true, colorlinks]{hyperref}
\hypersetup{pdffitwindow=true,linkcolor=blue,citecolor=blue,urlcolor=blue,menucolor=blue}

\begin{document} 

\rightline{\Huge \bf WHAT IS ... a Thin Group?}
\vskip1cm
\rightline{\it \huge Alex Kontorovich, D. Darren Long,} 
\vskip.1cm
\rightline{\it \huge Alexander Lubotzky, and Alan W. Reid%
\footnote{Alex Kontorovich is a professor of mathematics at Rutgers University. His email address is
{\tt alex.kontorovich@rutgers.edu}.
D. Darren Long is a professor of mathematics at the University of California, Santa Barbara. His email address is {\tt long@math.ucsb.edu}.
Alexander Lubotzky is the Maurice and Clara Weil professor mathematics at the Hebrew University of Jerusalem. His email address is {\tt alex.lubotzky@mail.huji.ac.il}.
Alan W. Reid is the Edgar Odell Lovett  professor of mathematics at Rice University. His email address is {\tt alan.reid@rice.edu}. 
The authors would like to gratefully acknowledge
an NSF FRG grant which facilitated their collaboration.
}}
\vskip1in






\begin{multicols}{2}


The group $\SL_2(\Z)$ of  $2\times 2$ integer matrices with unit determinant is a quintessential {\bf arithmetic group}. By this we mean that there is an  
{\bf algebraic group}, that is, a variety defined by polynomial equations, namely,
$$
\SL_2 \ : \ \{(a,b,c,d):ad-bc-1=0\},
$$
whose points over a ring
happen to also form a group (under standard matrix multiplication, which is a polynomial map in the entries); then $\SL_2(\Z)$ is the set of integer points in this algebraic group.
More generally, an arithmetic group $\G$ is a finite-index subgroup of the integer points $G(\Z)$ of an algebraic group $G$.
Roughly speaking, a ``thin'' group is an {\it infinite}-index subgroup of an arithmetic group which ``lives'' in the same algebraic group, as explained below.

While the term ``thin group''\footnote{Sometimes ``thin matrix group,'' not to be confused with other notions of thin groups in the literature.} was coined  in the last 10-15 years by Peter Sarnak, such groups had been studied as long as 100-150 years ago;
indeed,  they appear 
naturally
in the theory of Fuchsian and Kleinian groups.
For a long while, they
 were
 largely
  discarded as ``irrelevant'' to arithmetic, in part because there was not much one could do with them. 
More recently, thin groups have become a ``hot topic'' thanks to the explosion of activity in ``Super Approximation'' (see below). Armed with this new and massive hammer,
lots of previously unrelated problems 
in number theory, geometry, and group theory
started looking 
like nails. Our goal here is to describe some of these nails at a basic level; for a more advanced treatment of similar topics, the reader would do well to consult \cite{Sarnak2014}.

\medskip
Let's get to the general definition from first seeing some (non-)examples.
Take your favorite pair, $A, B$, say, of $2\times2$ matrices in $\SL_2(\Z)$ and let $\G=\<A,B\>$ be the group generated by them; should $\G$ be called thin?

\begin{Ex}\label{ex:1} \emph{
Suppose you choose $A=\mattwos1101$ and $B=\mattwos01{-1}0$. Then, as is well-known, $\G$ is all of $\SL_2(\Z)$. This cannot be  called  ``thin'';  it's the whole group.
}
\end{Ex}

\begin{Ex} \emph{
 If you choose $A=\mattwos1201$ and $B=\mattwos1021$, 
then the resulting $\G$  is also well-known to be a {\bf congruence group}, meaning roughly that the group is defined by congruence relations. More concretely, $\G$ turns out to be  the subset of $\SL_2(\Z)$ of all matrices with  diagonal entries congruent  to $1(\mod 4)$ and evens off the diagonal; it is a good exercise to check that these congruence restrictions do indeed form a group. 
It is not hard to prove that the index
\footnote{If the reader was expecting this index to be 6, that would be correct in $\PSL_2(\Z)$, or alternatively, if we added $-I$ to $\G$.} of  $\G$ in $\SL_2(\Z)$ is 12,
 so just 12 cosets of $\G$ will be enough to cover all of $\SL_2(\Z)$; that also doesn't qualify as  thin.
}
\end{Ex}

\begin{Ex}\label{ex:3} \emph{
 Say you chose $A=\mattwos1401$ and $B=\mattwos1601$;
that will generate $\G=\mattwos1{2\Z}01$, the group of  upper triangular matrices with an even upper-right entry.
This group is certainly of infinite index in $\SL_2(\Z)$, so now is it thin? Still no. The reason is that $\G$ fails to ``fill out'' the algebraic variety $\SL_2$. That is, there are ``extra'' polynomial equations satisfied by $\G$ besides $\det=1$; namely, $\G$ lives in the strictly smaller {\bf unipotent} (all eigenvalues are 1) algebraic group
$$
U  :  \{(a,b,c,d) : ad-bc-1=a-1= d-1=c=0\}.
$$
The fancy way of saying this is that $U$ is the {\bf Zariski-closure} of $\G$, written 
$$
U\ = \ \Zcl(\G).
$$ 
That is, $\Zcl(\G)$ is the algebraic group given by all polynomial equations satisfied by all elements of $\G$. And if we look at the integer points of $U$, we get $U(\Z)=\mattwos1\Z01$, in which $\G$ has finite index (namely, two). 
So again $\G$ is not thin.
}
\end{Ex}

\begin{Ex} \emph{
Take $A=\mattwos2111$ and $B=\mattwos5332$. This example is a little more subtle. The astute observer will notice that $B=A^2$, so $\G=\<A\>$, and moreover that 
$$
A^n=\mattwos{f_{2n+1}}{f_{2n}}{f_{2n}}{f_{2n-1}},
$$ 
where $f_n$ is the $n$-th Fibonacci number, determined by $f_{n+1}={f_{n}}+{f_{n-1}}$ and initialized by $f_0=f_1=1$.
Again it is easy to see that $\G$ is an infinite index subgroup of $\SL_2(\Z)$, and unlike \exref{ex:3}, all the entries are changing. But it is still not thin! 
Let $\phi={1+\sqrt{5}\over 2}$ be the golden mean and $K=\Q(\phi)$; there is a matrix $g\in\SL_2(K)$ which conjugates $\G$ to 
$$
\G_1\ = \ g\G g^{-1}\ = \ \{\mattwos{\phi^{2n}}00{\phi^{-2n}} : n\in\Z\}.
$$
The latter group lives inside the ``diagonal'' algebraic group
$$
D \ :  \ \{(a,b,c,d): b=c= ad-1=0\}.
$$
This group is an example of what's called an {\bf algebraic torus}; the group of complex points $D(\C)$ is isomorphic to the multiplicative ``torus'' $\C^\times$.
When we conjugate back, the torus $D$ goes to
$$
D_1\ = \ g^{-1}D\, g\  =\ \Zcl(\G);
$$ 
the variety $D_1$ is now defined by equations with coefficients in $K$, not $\Q$.
The rational integer points of $D_1$ are exactly $\G=D_1(\Z)$, so $\G$ is not a thin group.
}
\end{Ex}

\begin{Ex}\label{ex:5} \emph{
 This time, let $A=\mattwos1401$ and $B=\mattwos01{-1}0$, with $\G=\<A,B\>$. 
If we replace the upper-right entry $4$ in $A$ by $1$, we're back to \exref{ex:1}. So at first glance, perhaps this $\G$ has index $4$ or maybe $8$ in $\SL_2(\Z)$?
It turns out that $\G$ actually has infinite index (see, e.g., \cite[\S4]{Kontorovich2013} for a gentle discussion). What is its Zariski closure? Basically the only subvarieties of $\SL_2$ which are also groups look, up to conjugation, like $U$ and $D$ (and $UD
$), and it is easy to show  that $\G$ lives in no such group.
More generally, any subgroup of infinite index in $\SL_2(\Z)$ that is not virtually (that is, up to finite index) abelian is necessarily thin. Indeed, being non-virtually abelian rules out all possible proper sub-algebraic groups of $\SL_2$, implying that $\Zcl(\Gamma) = \SL_2$.
}
\end{Ex}


It is now a relatively simple matter to give an almost-general definition.

\begin{Def}
Let $\G<\GL_n(\Z)$ be a subgroup and let $G=\Zcl(\G)$ be its Zariski closure. We say $\G$ is a {\bf thin group} if the index of $\G$ in the integer points $G(\Z)$ is infinite. (Most people add that $\G$ should be finitely generated.)
\end{Def}

For more context, we return to the classical setting of a congruence group $\G<\SL_2(\Z)$. Such a group
acts on the upper half plane $\bH=\{z\in\C:\Im z>0\}$ by fractional linear transformations 
$$
\mattwo abcd:z\mapsto {az+b\over cz+d},
$$ 
and much 20th and 21st century mathematics has been devoted to the study of:
\begin{itemize}
\item ``Automorphic forms,'' meaning 
eigenfunctions 
$
\vf:\bH\to\C
$ 
of the hyperbolic Laplacian $\gD=y^2(\dd_{xx}+\dd_{yy})$ which are $\G$-automorphic, that is, 
$$
\vf(\g z)=\vf(z),
$$ 
for all $\g\in\G$ and $z\in\bH$,
and 
square-integrable (with respect to a certain invariant measure)
on the quotient $\G\bk\bH$. These are called ``Maass forms'' for Hans Maass's foundational papers in the 1940's. Their existence and abundance in the case of congruence groups is a consequence of the celebrated Selberg trace formula, developed in the 1950's.
\item ``$L$-functions'' attached to such $\vf$. These are certain ``Dirichlet series,'' meaning functions of the form 
$$
L_\vf(s) \ = \ \sum_{n\ge1}{a_\vf(n) \over n^s},
$$
where
 $a_\vf(n)$ is a sequence of complex numbers 
 called the
 ``Fourier coefficients'' 
 of $\vf$. When $\vf$ is also an eigenfunction of so-called ``Hecke operators'' and normalizing $a_\vf(1)=1$, these $L$-functions are also multiplicative, enjoying Euler products of the form
$$
L_\vf(s) \ = \ \prod_p\left(1+{a_\vf(p) \over p^s}+{a_\vf(p^2) \over p^{2s}}+\cdots\right)
,
$$
where the product runs over primes. Needless to say, such $L$-functions are essential in modern analytic number theory, with lots of fascinating applications to primes and beyond.
\item
More generally, one can define related objects (called ``automorphic representations'')  on other arithmetic groups $G(\Z)$, and study their $L$-functions. The transformative insight of Langlands is the conjectured interrelation of these  on different groups, seen most efficiently through the study of operations on their $L$-functions. Consequences of these hypothesized interrelations include the Generalized Ramanujan and Sato-Tate Conjectures, among many, many others. We obviously have insufficient capacity to do more than graze the surface here.
\end{itemize}

Hecke, in studying \exref{ex:5}, found that his theory of Hecke operators fails for thin groups, %
 so such $L$-functions would not have Euler products,\footnote{Without Euler products, $L$-functions can have zeros in the region of absolute convergence; that is, the corresponding Riemann Hypothesis can fail dramatically!} and hence no direct applications to questions about primes. Worse yet, the Selberg trace formula breaks down, and  there are basically no Maass forms to speak of (nevermind the $L$-functions!). So for a long while, it seemed like  thin groups, although abundant, did not appear particularly relevant to arithmetic problems.
 
About 15 years ago, a series of stunning breakthroughs led to the theory of ``Super Approximation,'' as described below, 
and for the first time allowed a certain Diophantine analysis on thin groups, from which many striking applications soon followed. To discuss these, we first describe the more classical theory of
 {\em Strong Approximation}. In very rough terms, this theory says that from a certain algebraic perspective, ``thin groups are indistinguishable from their arithmetic cousins,'' by which we mean the following.
 
It is not hard\footnote{Though if you think it's  completely trivial, try finding a matrix $\g\in\SL_2(\Z)$ whose reduction mod $5$ is, say, $\mattwos2013$. The latter is indeed an element of $\SL_2(\Z/5\Z)$, since it has determinant $6\equiv1(\mod 5)$.}  
to show that reducing $\SL_2(\Z)$ modulo a prime $p$ gives all of $\SL_2(\Z/p\Z)$. 
What happens if we reduce the group $\G$ in \exref{ex:5} mod $p$? Well, for $p=2$, we clearly have a problem, since the generator $A=\mattwos1401$ collapses to the identity. But for any other prime $p\neq2$, the integer $4$ is a unit (that is, invertible mod $p$), so some power of $A$ is congruent to $\mattwos1101$ mod $p$. Hence on reduction mod  (almost)  {\it any} prime,  we cannot distinguish \exref{ex:5} from \exref{ex:1}!
That is, even though $\G$ in \exref{ex:5} is thin, the reduction map $\G\to\SL_2(\Z/p\Z)$ is onto. 
The Strong Approximation theorem \cite{MatthewsVasersteinWeisfeiler1984} says that if $\G<\SL_n(\Z)$ has, say, full Zariski closure $\Zcl(\G)=\SL_n$,  then $\G\to\SL_n(\Z/p\Z)$ is onto for all but finitely many primes $p$. 
In fact, this reasoning can be reversed, giving a very easy check of Zariski density: 
 if for a single prime $p\ge5$, the reduction of $\G$ mod $p$ is all of $\SL_n(\Z/p\Z)$, then the Zariski closure of $\G$ is automatically all of $\SL_n$; see \cite{Lubotzky1999} for details.

One immediate
 caveat 
 is that, if one is not careful,
 Strong Approximation can fail. For a simple example,
 try finding a $\g\in\GL_2(\Z)$ which mod $5$ gives $\mattwos1234$.
The problem is that 
   $\GL_{n}(\Z)$ does not map onto $\GL_{n}(\Z/p\Z)$,
 since the only determinants of the former are $\pm1$, while the latter has determinants in all of $(\Z/p\Z)^{\times}$. But such 
 obstructions are well-understood and classical. (In fancy language, $\GL_n$ is reductive, while $\SL_n$ is semisimple.)


For Super Approximation, we study not only {\it whether} these generators $A$ and $B$ in \exref{ex:5} fill out $\SL_2(\Z/p\Z)$, but the more refined question of {\it how rapidly} they do so.
To quantify this question, construct for each (sufficiently large) prime $p$ the {\it Cayley graph}, $\sG_p$, whose vertices are the elements of $\SL_2(\Z/p\Z)$ and two vertices (i.e., matrices) are connected if one is sent to the other under one of  the four generators $A^{\pm1}$, $B^{\pm1}$. When $p=3$, the graph is:%
\footnote{This graph is begging us to identify each node $\g$ with $-\g(\mod p)$, that is, work in $\PSL_2$.}
\begin{center}
\includegraphics[width=.4\textwidth]{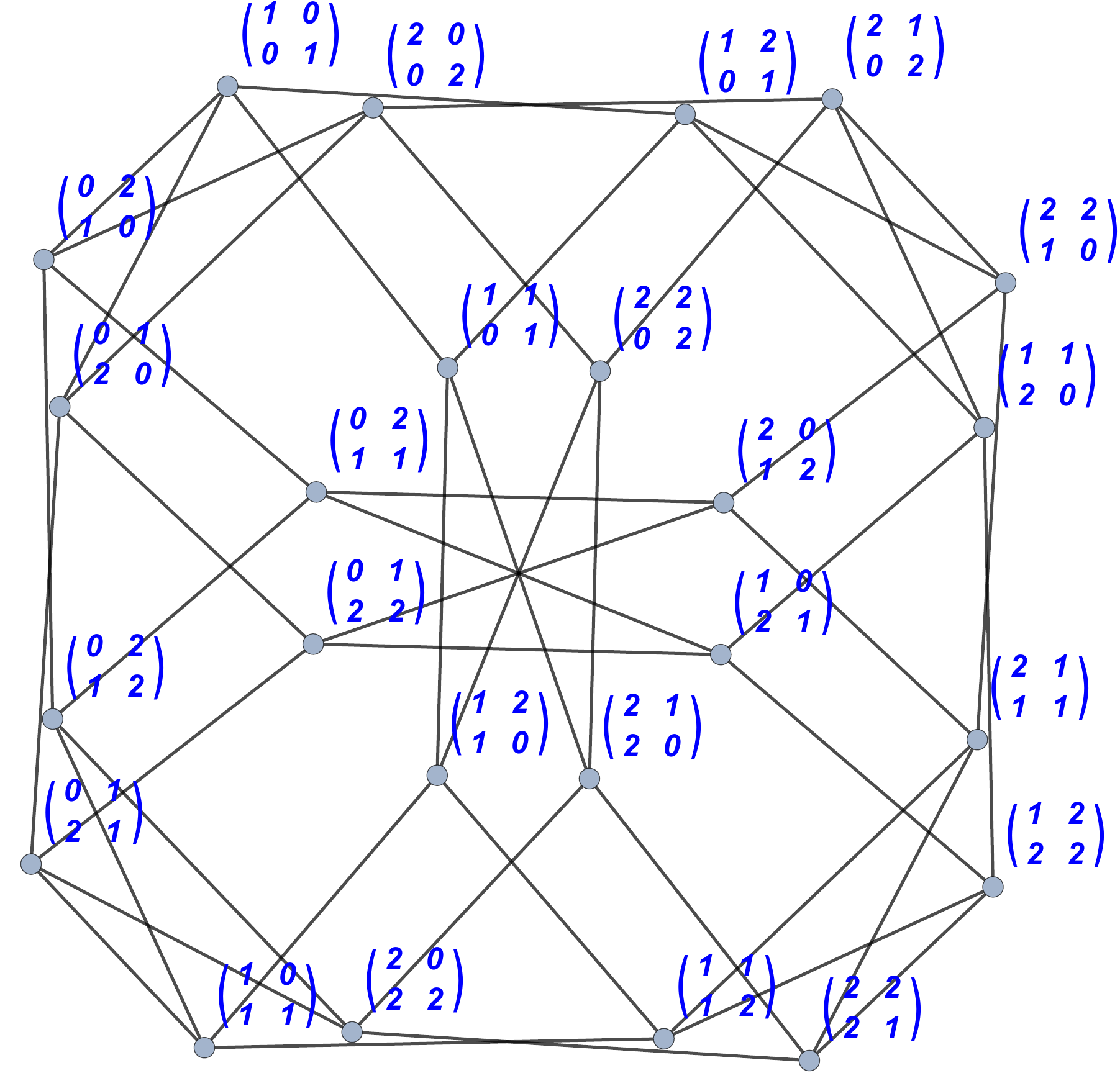}
\end{center}
This is a $k$-regular graph with $k=4$, that is, every vertex $\g\in\SL_2(\Z/p\Z)$ has four neighbors.
The ``graph Laplacian''  of $\sG_p$ is the matrix
$
\gD\ : = \ I-\frac1k \cA,
$ 
where $\cA$ is the adjacency matrix of the graph. By the  {\it spectrum} of $\sG_p$, we mean  the eigenvalues
$$
\gl_0^{(p)}\le \gl_1^{(p)}\le\cdots
$$
of $\gD$.
In the case of the graph 
above,
the spectrum is:
$$
\left\{0,\frac{1}{2},\frac{1}{2},\frac{1}{2},\frac{1}{2},\frac{3}{4},\frac{3}{4},\frac{5}{4},\cdots, \frac{1}{8} \left(7+\sqrt{17}\right)\approx1.39\right\}
$$
Notice that the bottom eigenvalue $\gl_0$ is $0$ (corresponding to the constant function), and has multiplicity $1$; this is due to Strong Approximation -- the graph is connected!  (In general, the multiplicity of the bottom eigenvalue is the number of connected components.)
Hence the first eigenvalue above the bottom, $\gl_1^{(p)}$, is strictly positive, which by standard techniques implies that a random walk on the graph is ``rapidly mixing'' (see, e.g., \cite{DavidoffSarnakValette2003}).
But we have infinitely many graphs $\sG_p$, one for each prime, and  {\it a priori}, it might be the case that the mixing rate deteriorates as $p$ increases.
Indeed, $\gl_1^{(p)}$ goes from $\foh$ when $p=3$ down to $\gl_1^{(p)}\approx0.038$ when $p=23$, for which the graph has about $12,000$ vertices.
Super Approximation is precisely the statement that this deterioration does not continue indefinitely: there exists some $\vep>0$ so that, for all $p$,
$$
\gl_1^{(p)}\ge \vep.
$$
That is, the rate of  mixing is {\it uniform} over the entire {\it family} of Cayley graphs $\sG_p$. 
(This is what's called an {\it expander} family, see \cite{Sarnak2004, Lubotzky2012}.)

For congruence groups, Super Approximation
is now a classical fact: it is a consequence of ``Kazhdan's Property T'' in higher rank (e.g. for groups like $\SL_n(\Z)$ $n\geq 3$), and of  non-trivial bounds towards the ``Generalized Ramanujan Conjectures'' in rank one (for example,
isometry groups of hyperbolic spaces); see, e.g.,  \cite{Lubotzky1994, Sarnak2005} for an exposition.
A version of Super Approximation for some more general (arithmetic but not necessarily congruence) lattices was established by Sarnak-Xue \cite{SarnakXue1991}.

For thin subgroups $\G<\SL_{n}(\Z)$, major progress was made 
by Bourgain-Gamburd \cite{BourgainGamburd2008},
who established Super Approximation (as formulated above) for  $\SL_{2}$.
This built on
a sequence of striking results in  Additive Combinatorics, namely the Sum-Product Theorem \cite{BourgainKatzTao2004} and Helfgott's Triple Product Theorem \cite{Helfgott2008},
and
 prompted a slew of activity by many people (e.g. \cite{Varju2012, PyberSzabo2016, BreuillardGreenTao2011}), culminating in an (almost) general Super Approximation theorem of Salehi-Golsefidy and Varju \cite{SalehiVarju2012}.

Simultaneously, it was realized that many natural problems in number theory, groups, and geometry {\it require} one to treat these aspects of thin (as opposed to arithmetic) groups. 
Two quintessential such, discussed at length in \cite{Kontorovich2013}, are the Local-Global Problem for integral Apollonian packings \cite{BourgainKontorovich2014a}, and Zaremba's conjecture on ``badly approximable'' rational numbers \cite{BourgainKontorovich2014}.
Other related problems subsequently connected to thin groups (see the exposition in \cite{Kontorovich2016}) include McMullen's Arithmetic Chaos Conjecture and a problem of Einsiedler-Lindenstrauss-Michel-Venkatesh on low-lying fundamental geodesics on the modular surface.
The latter problem, eventually resolved  in \cite{BourgainKontorovich2017},
was the catalyst for the development of the Affine Sieve \cite{BourgainGamburdSarnak2010, SalehiSarnak2013}; see more discussion in \cite{Kontorovich2014}.
Yet a further direction was opened by the realization that the Affine Sieve can be extended to
what may be
called the  ``Group Sieve,''
used to great effect on problems in Group Theory and Geometry in, e.g., \cite{Rivin2008, Kowalski2008, LongLubotzkyReid2008, LubotzkyMeiri2012}.
We will not rehash these topics, 
choosing instead to
end by
highlighting the difficulty of answering the slight rewording of the title: 
$$
\text{\it
Can you tell... whether a given group is Thin?
}
$$

\begin{Ex}\label{ex:7}
\emph{
To ease us into a higher rank example, consider the group $\G<\SL_3(\Z)$ generated by 
$$
A
=
\bp
1&1&0\\
0&1&0\\
0&0&1
\ep
\text{
and
}
B
=
\bp
0&1&0\\
-1&0&0\\
0&0&1
\ep.
$$
A moment's inspection reveals that $\G$ is just a copy of $\SL_2(\Z)$ (see \exref{ex:1}) in the upper left $2\times2$ block of $\SL_3$. This $\G$ has Zariski closure isomorphic to $\SL_2$, and is hence not thin.
}
\end{Ex}

\begin{Ex}\label{ex:8}\emph{
Here's a much more subtle example. Set
$$ A = \left( \begin{array}{ccc}0  &  0 & 1\\                                                           
1 &0 &0 \\                                                                                                    
0 &1 &0 \end{array}\right) $$
and
$$
B=\left( \begin{array}{ccc}1  &   2  & 4\\
0 & - 1 & -1  \\
0 & 1  & 0 \end{array}\right). 
$$
It is not hard to show that group $\G=\<A,B\>$ 
has full Zariski closure, $\Zcl(\G)=\SL_3$. Much more striking (see \cite{LongReidThistlethwaite2011}) is that $\G$ is 
  a {\em faithful} representation of the  ``$(3,3,4)$ hyperbolic triangle'' group 
  $$
  T= \<A,B:A^3=B^3=(AB)^4=1\>
  $$ 
  into $\SL_3(\Z)$;  that is, the generators have these relations and no others. It then follows that  $\G$ is necessarily of infinite index in $\SL_3(\Z)$, that is, thin.
}\end{Ex}

\begin{Ex}\label{ex:9}\emph{
The 
matrices
$$
A=
\bp
0&0&0&-1\\
1&0&0&-1\\
0&1&0&-1\\
0&0&1&-1
\ep
,\ 
B=
\bp
1&0&0&5\\
0&1&0&-5\\
0&0&1&5\\
0&0&0&1
\ep
$$
generate a group $\G<\SL_4(\Z)$ whose Zariski closure turns out to be the symplectic group $\Sp(4)$. 
The interest in these particular matrices is that they generate the ``monodromy group'' of a certain (Dwork) hypergeometric equation. 
It was shown in \cite{BravThomas2014} that this group is thin. For general monodromy groups, determining who is thin or not is wide open; see related work in \cite{Venkataramana2014} and \cite{FuchsMeiriSarnak2014}, as well as the discussion in \cite[\S3.5]{Sarnak2014}.
}\end{Ex}

\begin{Ex}\label{ex:10}\emph{
The four
matrices
$$
\left(
\begin{array}{cccc}
 1 & 0 & 0 & 0 \\
 0 & 1 & 0 & 0 \\
 0 & 0 & -1 & 0 \\
 0 & 0 & 0 & 1 \\
\end{array}
\right),
\left(
\begin{array}{cccc}
 1 & 0 & 0 & 0 \\
 1 & 1 & 1 & 0 \\
 -2 & 0 & -1 & 0 \\
 0 & 0 & 0 & 1 \\
\end{array}
\right),
$$
$$
\left(
\begin{array}{cccc}
 0 & 1 & 0 & 0 \\
 1 & 0 & 0 & 0 \\
 0 & 0 & 1 & 0 \\
 0 & 0 & 0 & 1 \\
\end{array}
\right),
\left(
\begin{array}{cccc}
 3 & 2 & 0 & 1 \\
 2 & 3 & 0 & 1 \\
 0 & 0 & 1 & 0 \\
 -12 & -12 & 0 & -5 \\
\end{array}
\right)
$$
generate a group $\G<\GL_4(\Z)$.
Its Zariski closure turns out to be the ``automorphism group'' 
of a certain quadratic form of signature $(3,1)$.
By a standard process (see, e.g., \cite[p. 210]{Kontorovich2013}), such a $\G$  acts on hyperbolic 3-space 
$$
\bH^3=\{(x_1,x_2,y):x_j\in\R,y>0\},
$$ 
and in this action, each matrix represents inversion in a hemisphere. 
These inversions are shown in red in the figure below, as is the set of  limit points of a $\G$ orbit (viewed in the boundary plane $\R^2=\{(x_1,x_2,0)\}$); the latter turns out to be a fractal  circle packing:
\begin{center}
\includegraphics[width=.45\textwidth]{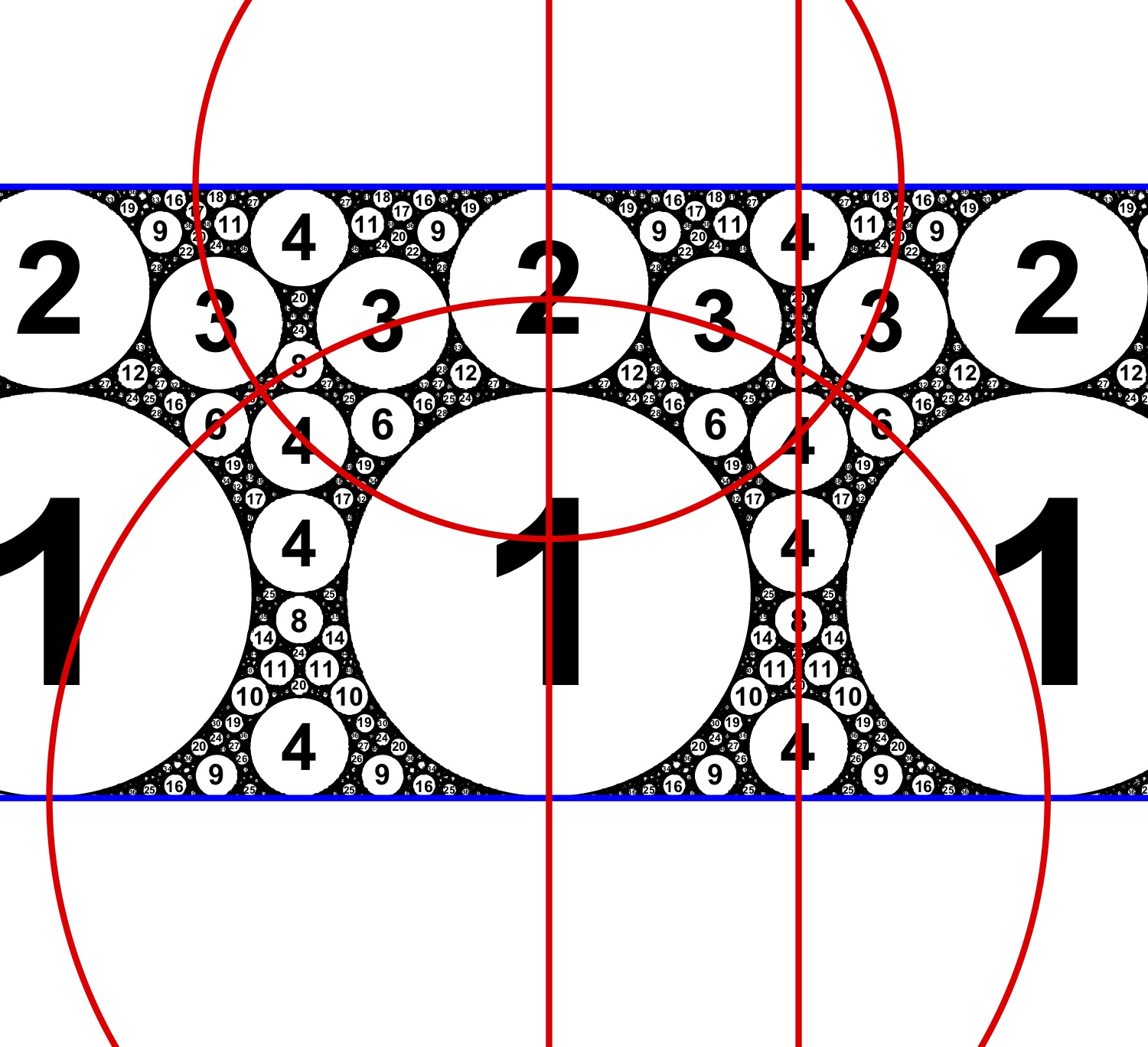}
\end{center}
Here circles are labeled with the reciprocal of their radii (notice these are all integers!).
This limit set is an example of a ``crystallographic packing,'' introduced (and partially classified) in
\cite{KontorovichNakamura2018} as a vast generalization of integral Apollonian circle packings. It follows from the fractal nature of this limit set that $\G$ is indeed a thin group.
}\end{Ex}

In a sense that can be made precise (see \cite{LubotzkyMeiri2012, Aoun2011, Rivin2010, FuchsRivin2017}),
random subgroups of arithmetic  groups are thin. But 
lest we leave the reader with the false impression that the theory is truly well-developed and on solid ground, we demonstrate our ignorance with the following basic challenge.

\begin{Ex} 
\emph{
The following group arises naturally through certain geometric considerations in \cite{LongReidThistlethwaite2011}:
let $\G=\<A,B\><\SL_3(\Z)$ with
$$
A = \left( \begin{array}{ccc}1  & 1 &  2 \\ 0 & 1 &  1 \\ 0 & -3 & -2 \end{array}\right),
\
B = \left( \begin{array}{ccc} -2 &  0 &  -1 \\ -5 & 1 &  -1 \\ 3 & 0& 1 \end{array}\right).$$
Reduced mod 7, this $\G$ is all of $\SL_3(\Z/7\Z)$, so its Zariski closure is $\SL_3$. Is it thin? 
As of this writing, nobody knows! 
}
\end{Ex}

\end{multicols}
\bibliographystyle{alpha}
\bibliography{AKbibliog}
\end{document}